\documentclass{gen-j-l}%
\usepackage{amsfonts}
\usepackage{geometry}%
\usepackage{amsmath}%
\usepackage{amssymb}%
\usepackage{graphicx}
\newtheorem{theorem}{Theorem}[section]

\theoremstyle{definition}
\newtheorem{definition}[theorem]{Definition}

\theoremstyle{remark}

\numberwithin{equation}{section}
\theoremstyle{plain}

\copyrightinfo{2017}{Carey Caginalp}
\begin{document}
\title[Conservation Laws With Random and Deterministic Data]{Conservation Laws With Random and Deterministic Data}
\author{Carey Caginalp}
\address{182 George St, Box F, Division of Applied Mathematics, Brown University 02912}
\thanks{The author thanks Professors Menon and Dafermos and Dr. Kaspar for valuable discussions. This work was partially supported by NSF grants DMS 1411278 and DMS 1148284 as well as the NSF Graduate Research Fellowship.}
\email{carey\_caginalp@brown.edu}
\urladdr{http://www.pitt.edu/\symbol{126}careycag/}

\begin{abstract}
The dynamics of nonlinear conservation laws have long posed fascinating
problems. With the introduction of some nonlinearity, e.g. Burgers' equation,
discontinuous behavior in the solutions is exhibited, even for smooth initial
data. The introduction of randomness in any of several forms into the initial
condition makes the problem even more interesting. We present a broad spectrum
of results from a number of works, both deterministic and random, to provide a
diverse introduction to some of the methods of analysis for conservation laws.
Some of the deep theorems are applied to discrete examples and illuminated
using diagrams.

\end{abstract}
\maketitle

\section{Introduction}

\subsection{\textbf{Background.}}

In the effort to create a mathematical description of turbulence, an important
building block is the study of shocks and rarefactions together with random
initial conditions. Although the model is a somewhat coarse description of
turbulence in practice, Burgers' equation is extensively studied \cite{HP,
LA1, LA2} as a test case for new methods and types of randomness. It also
possesses the surprising feature of producing discontinuous solutions, even
from smooth initial data. From there it is then reasonable to seek broader
classes of equations to which these properties can be extended.

The link between many-particle systems and fluid mechanics, shock waves, and
PDEs poses important problems in understanding the continuum limit. In this
paper, we provide a summary of various results in the field and potential
directions for open problems in the future. Our aim is to tie together a
number of vastly different approaches across the scope of kinetic theory
involving conservation laws (not limited simply to the widely studied Burgers'
equation) in a comprehensive note that serves both as a general introduction
and as a starting point for readers who may wish to delve into the more
technical aspects in the references herein. We also highlight the role of
various types of randomness in the initial conditions, and the preservation
(or lack thereof) of certain prescribed structure in the solution as time
advances. We illuminate the theory with a series of discrete examples.
Randomness in these initial conditions, together with formation and
interaction of resulting shocks, forms a basic model that is a first step
toward a mathematical description of turbulence. This is subsequently useful
in a wide array of applications in fields such as engineering, one such
application being attempts to control turbulent flows.

\subsection{\textbf{Burgers' Equation Derived From Pressureless Limit.}}

In the continuum, one represents gas dynamics in one dimension by density and
velocity fields, $\rho\left(  t,x\right)  $ and $u\left(  t,x\right)  ,$
respectively. \ The Euler equations, given by%
\begin{align}
\rho_{t}+\left(  \rho v\right)  _{x}  & =0\text{ (conservation of
mass),}\nonumber\\
\left(  \rho v\right)  _{t}+\left(  \rho v^{2}+p\right)  _{x}  & =0\text{
(conservation of momentum),}\nonumber\\
\left(  \rho E\right)  _{t}+\left(  \rho Ev+pv\right)  _{x}  & =0\text{
(conservation of energy),}\label{Eulereq}%
\end{align}
provide a starting point for the approach of Brenier and Grenier \cite{BG},
which we detail in Section 2. One then considers inelastic collisions (under
which kinetic energy is not conserved) and take the pressureless limit,
formally obtaining the system%
\begin{align}
\partial_{t}\rho+\partial_{x}\left(  \rho u\right)   & =0\nonumber\\
\partial_{t}\left(  \rho u\right)  +\partial_{x}\left(  \rho u^{2}\right)   &
=0\label{system}%
\end{align}
In working with this system, Radon measures (see \cite{RY}, p. 455) provide a
key tool in making the interpretation of the equations precise in the most
general case. \ The appropriate system of conservation laws provides
conservation of mass and momentum, and the unknowns include these quantities
for specific particles, along with velocity. \ Under some basic assumptions,
the mathematical tool of the Radon-Nikodym derivative, which is essentially
the derivative of a measure \cite{RY} (p. 385), then allows them to define in
a rigourous sense the quotient of these measures (momentum and mass), a
mathematical analog to velocity being the quotient of momentum over mass.

\subsection{\textbf{Analyzing Burgers' Equation From Another Perspective: Flow
Maps.}}

In Section 3, we analyze a similar problem with a very different approach. The
same system from \cite{BG} is presented in E, Rykov, and Sinai \cite{ERS}, and
is equivalent, under smooth solutions, to Burger's equation and a transport
equation%
\begin{align}
u_{t}+\left(  \frac{u^{2}}{2}\right)  _{x}  & =0\nonumber\\
\rho_{t}+\left(  \rho u\right)  _{x}  & =0
\end{align}
From here one can define%
\begin{equation}
x=\varphi_{t}\left(  y\right)  =y+tu_{0}\left(  y\right)  .
\end{equation}
Formal inversion of this flow map yields%
\begin{equation}
u\left(  x,t\right)  =u_{0}\left(  \varphi_{t}^{-1}\left(  x\right)  \right)
,\rho\left(  x,t\right)  =\rho_{0}\left(  \varphi_{t}^{-1}\left(  x\right)
\right)  \left\vert \frac{\partial x}{\partial y}\right\vert ^{-1}.
\end{equation}
The issue is that this flow map is only invertible up to some time $t$, at
which point $\varphi_{t}\left(  x\right)  $ is no longer one-to-one and the
inverse not well-defined. \ Specifically, we have a whole interval mapped into
one point, forming a shock. \ However, we note that the map $\varphi_{t}$
still defines a partition of the real line (mapping some points into points
and sometimes intervals into a single point). \ It is from this observation
that E, Rykov, and Sinai construct solution formulas (in the appropriate
sense, to be defined later) given by%
\begin{equation}
\varphi_{t}\left(  y\right)  =\frac{\int_{C_{t}\left(  y\right)  }\left(
\eta+tu_{0}\left(  \eta\right)  \right)  dP_{0}\left(  \eta\right)  }%
{\int_{C_{t}\left(  y\right)  }dP_{0}\left(  \eta\right)  },u\left(
x,t\right)  =\frac{\int_{D_{t}\left(  x\right)  }u_{0}\left(  \eta\right)
dP_{0}\left(  \eta\right)  }{\int_{D_{t}\left(  x\right)  }dP_{0}\left(
\eta\right)  }%
\end{equation}
where $D_{t}\left(  x\right)  ,C_{t}\left(  x\right)  $ have to do with the
aforementioned partition.

\subsection{\textbf{Entropy Solution and Variational Approach Using Stieltjes
Integral.}}

A similar variational approach is studied in Huang and Wang \cite{HW}, where
an entropy solution is constructed from a generalised potential%
\begin{equation}
F\left(  y;x,t\right)  =\int_{0+0}^{y-0}tu_{0}\left(  \eta\right)
+\eta-xdm_{0}\left(  \eta\right)  .
\end{equation}
They construct a set $S\left(  x,t\right)  $ to serve the role of determining
what point or interval is mapped into a point $\left(  x,t\right)  $. \ Here
the notation denotes a special kind of integral known as the \textit{Stieltjes
integral}, whereby integration is from the right limit of $0$ to the left
limit of $y$. Unlike the traditional Riemann or Lebesgue integrals, this may
result in different values for the expression, as shown in Section 4.

\subsection{\textbf{Introduction of Random Initial Conditions for More General
Conservation Laws.}}

In the work of Menon and Srinivasan \cite{MS} and Kaspar and Rezakhanlou
\cite{KR}, further analysis on these equations was performed to gain deeper
understanding of particle dynamics and the evolution of the system starting
with random initial data. The analysis went beyond Burgers' equation to the
more general case of a $C^{1},$ convex flux, and considered initial data that
was random rather than deterministic. Here, the initial conditions were
restricted to processes that were \textit{spectrally negative}, that is,
stochastic processes with jumps but only in one direction. In particular, only
downward jumps were permitted, as upward jumps lead to an immediate breakdown
of the statistics, in that other behavior such as rarefaction is observed.
Furthermore, this initial stochastic process was also assumed to be Markov.
Using the Levy-Khinchine representation for the Laplace exponent and other
methods, formal calculations were used to establish a number of results. One
remarkable assertion proven states that if one starts with strong Markov,
spectrally negative initial data, then this Markov property persists in the
entropy solution for any positive time $t>0$. This closure result can also
hold for a non-stationary process and can obtain an equation for a generator
in the case of specific kinds of stochastic processes. These results are
described in Section 5.

In a few special cases, such as Burgers' equation under white noise initial
conditions in \cite{FM}, one has the remarkable achievement of a closed form
solution up to the level of special functions. By starting from the full
Burgers' equation and taking the vanishing viscosity limit, one is lead to the
variational blueprint for this set of exact results, along with a crisp
geometric visualisation. We will examine this case in further detail in
Section 6.

\subsection{\textbf{Extension of Results to Flux Functions With Lesser
Regularity.}}

Subsequently, we extend these results to a further class of nonlinear flux
functions, providing results and derivations of two different hierarchies.
These were checked rigourously against various examples involving Riemann
initial data, a protoypical but important base case that is the building block
to more complicated or even random initial conditions. In particular, the
second hierarchy derived shows consistency through shock interactions without
any sort of extraneous resetting or additional conditions imposed, a feature
absent in many classical methods. This is presented in Section 7.

\subsection{\textbf{Open Problems.}}

Finally, in Section 8, we discuss future work in these areas, including the
possibility of proving a rigourous closure theorem for these hierarchies.
Other applications such as more computational testing of these equations under
various forms of random initial data are also open problems.

\section{Burgers' Equation And the Sticky Particle Model}

Brenier and Grenier \cite{BG} consider Burgers' equation applied to a
pressureless gas, described at a discrete level by a large collection of
sticky particles. The "sticky" part of the description corresponds to the
particles remaining together after inelastic collisions in accordance with
conservation of mass and momentum (but notably, \textit{not }conservation of
energy). At a continuum level, the model is described by density and velocity
fields $\rho\left(  t,x\right)  $ and $u\left(  t,x\right)  ,$ respectively,
that must satisfy%
\begin{align}
\partial_{t}\rho+\partial_{x}\left(  \rho u\right)   & =0\label{contsp}\\
\partial_{t}\rho u+\partial_{x}\left(  \rho u^{2}\right)   & =0.\nonumber
\end{align}
The system (\ref{contsp}) follows from formally letting the pressure go to
zero in (\ref{Eulereq}) or letting the temperature approach zero in the
Boltzmann equation. In particular, the system (\ref{contsp}) can be shown to
be equivalent to the inviscid Burgers' equation under the assumptions of
smooth solutions and positive densities.

However, for the considerations of the sticky particle model, such a reduction
is not as immediate. Application of (\ref{contsp}) to this problem presents
several difficulties:\ (i) under discrete particles, the fields are no longer
functions but must be considered as measures, (ii) the velocity field needs to
be well-defined almost everywhere with respect to the measure prescribed by
$\rho$, and (iii) the system must be supplemented by some entropy conditions.
In addition, an obvious choice such as the condition%
\begin{equation}
\partial_{t}\left(  \rho U\left(  u\right)  \right)  +\partial_{x}\left(  \rho
uU\left(  u\right)  \right)  \leq0
\end{equation}
for any smooth, convex $U$ is shown to be insufficient to guarantee uniqueness.

In this section we aim to illustrate how the flux function is constructed by
considering a discrete example with a finite number of particles. We also show
detail how the solution of the conservation law is linked to the
Hamilton-Jacobi equation in a viscosity sense.

\subsection{Discrete Example for Burgers' Equation}

\begin{figure}
[h]
\begin{center}
\includegraphics[width=\linewidth]{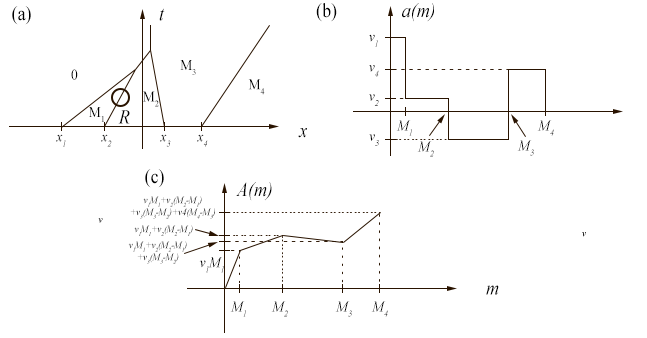}
\caption{(a) By taking a cross-section in time, one can obtain a cumulative
distribution function of the mass as a function of position; (b) Illustration
of the potential as a function of mass; (c) Illustration of the flux function
of mass.}
\label{BGfigure}
\end{center}
\end{figure}

For definiteness, take $n=4$ particles with masses $\left\{  m_{i}\right\}
_{i=1}^{4}$ at positions $\left\{  x_{i}\right\}  _{i=1}^{4}$ and with initial
velocities $\left\{  v_{i}\right\}  _{i=1}^{4},$, respectively, given by%
\begin{align}
m_{1}^{0}  & =\frac{1}{4},\text{ \ }m_{2}^{0}=\frac{1}{4},\text{ \ }m_{3}%
^{0}=\frac{1}{3},\text{ \ }m_{4}^{0}=\frac{1}{6},\nonumber\\
v_{1}^{0}  & =2,\text{ \ }v_{2}^{0}=1,\text{ \ }v_{3}^{0}=-\frac{1}{2},\text{
\ }v_{4}^{0}=1,\nonumber\\
x_{1}^{0}  & =-3,\text{ \ }x_{2}^{0}=-2,\text{ \ }x_{3}^{0}=1,\text{ \ }%
x_{4}^{0}=3,
\end{align}
These are plotted in Figure \ref{BGfigure}(a) along with characteristics as
the dynamics evolve in time. \ Our goal is to show how we build a PDE to model
the dynamics of this system, and that the solution of this PDE matches our
intuition about what should occur with sticky particles. \ Define for
notational convenience $M_{i}=\sum_{j=1}^{i}m_{i}$, $M_{0}=0$. \ The initial
distribution is given by $M^{0}\left(  x\right)  =\rho^{0}\left(
(-\infty,x]\right)  $ where $\rho^{0}$ jumps by $\Delta M_{i}$ at the points
$x_{i}$ (so clearly $M^{0}\left(  x_{i}-\right)  =M_{i-1}$, $M^{0}\left(
x_{i}+\right)  =M_{i}$). \ Now we define the function $a\left(  m\right)
=u_{0}\left(  x\right)  $ for $M^{0}\left(  x-\right)  \leq m\leq M^{0}\left(
x+\right)  .$ \ In this case, this simplifies to%
\begin{equation}
a\left(  m\right)  :=v_{i}\text{, }M_{i-1}\leq m\leq M_{i}\label{littleadef}%
\end{equation}
and is plotted in Fig. \ref{BGfigure}(b).

Now, we construct the flux function:%
\begin{equation}
A\left(  m\right)  =\int_{0}^{m}a\left(  m^{\prime}\right)  dm^{\prime}%
\end{equation}
and plot it in Figure \ref{BGfigure}(c).

\ 

We want to construct a weak solution to the differential equation%
\begin{equation}
\partial_{t}M+\partial_{x}\left(  A\left(  M\right)  \right)  =0
\end{equation}
and show that it works for the flux function proposed above. \ A weak solution
will satisfy%
\begin{equation}
\int_{-\infty}^{\infty}\int_{0}^{\infty}\varphi_{t}M+\varphi_{x}A\left(
M\right)  =0
\end{equation}
for every smooth function $\varphi\in C_{c}^{\infty}\left(  \mathbb{R\times
}\left(  0,\infty\right)  \right)  .$ \ We want our flux function to satisfy
the Rankine-Hugoniot condition at each shock. \ This means if we have
$M=M_{l}$ as the left-hand limit and $M=M_{r}$ as the right-hand limit,
$A\left(  M\right)  $ should satisfy%
\begin{equation}
\frac{A\left(  M_{r}\right)  -A\left(  M_{l}\right)  }{M_{r}-M_{l}}=\sigma
\end{equation}
where $\sigma$ is the slope of the parametrised curve describing the shock.
\ In particular, choosing a test function with compact support in the region
$R$ sketched in Fig. 1(a) yields%
\begin{equation}
\frac{A\left(  M_{r}\right)  -A\left(  M_{l}\right)  }{M_{r}-M_{l}}%
=\frac{v_{1}M_{1}+v_{2}\left(  M_{2}-M_{1}\right)  -v_{1}M_{1}}{M_{2}-M_{1}%
}=\frac{v_{2}m_{2}}{m_{2}}=v_{2}.
\end{equation}
Hence, the Rankine-Hugoniot condition is satisfied. \ Similarly, one verifies
it holds at the other discontinuities.

\subsection{Verifying that the potential $\Psi$ is a viscosity solution}

Another result of Brenier and Grenier involves linking this problem with the
Hamilton-Jacobi equation. \ To this end, they define a potential%
\begin{equation}
\Psi\left(  x,t\right)  =\int_{-\infty}^{x}M\left(  t,y\right)  dy.
\end{equation}
This is a viscosity solution in the sense of Crandall-Lions of the following
Hamilton-Jacobi equation:%
\begin{equation}
\partial_{t}\Psi+A\left(  \partial_{x}\Psi\right)  =0
\end{equation}
and is derived from the second Hopf formula. \ Indeed, in Bardi and Evans
[BE], the following is proven.

\begin{theorem}
Assume $A:\mathbb{R}^{n}\rightarrow\mathbb{R}$ is continuous and $\Psi
_{0}:\mathbb{R}^{n}\rightarrow\mathbb{R}$ is uniformly Lipschitz and convex.
\ Then%
\begin{equation}
\Psi\left(  x,t\right)  =\sup_{y}\inf_{z}\left\{  u_{0}\left(  z\right)
+y\cdot\left(  x-z\right)  -tA\left(  y\right)  \right\}  \text{ (for
}t>0\,,\text{ }x\in\mathbb{R}\text{)}%
\end{equation}
is the unique uniformly continuous viscosity solution of%
\begin{align}
\Psi_{t}+A\left(  \partial_{x}\Psi\right)   & =0\text{ in }\mathbb{R}%
^{n+1}\nonumber\\
\Psi\left(  \cdot,0\right)   & =\Psi_{0}\left(  \cdot\right)
\end{align}

\end{theorem}

This is valuable because it provides a solution to the Hamilton-Jacobi
equation without the usual convexity assumptions on the flux function $A$.
\ Performing some rearrangements, we have that%
\begin{equation}
\Psi\left(  x,t\right)  =\sup_{m}\left\{  \inf_{z}\left(  u_{0}\left(
z\right)  -m\cdot z\right)  +m\cdot x-tA\left(  y\right)  \right\}  =\sup
_{m}\left\{  -\Psi_{0}^{\ast}\left(  m\right)  +m\cdot x-tA\left(  m\right)
\right\}  ,
\end{equation}
where $\Psi_{0}^{\ast}$ is the Legendre-Fenchel transform given by%
\begin{equation}
\Psi_{0}^{\ast}\left(  m\right)  =\inf\left\{  u_{0}\left(  z\right)  -y\cdot
z\right\}
\end{equation}
i.e., in the notation of Brenier and Grenier \cite{BG}, equation (25)%
\begin{equation}
\Psi\left(  t,x\right)  =\sup_{0\leq m\leq1}\left\{  xm-\Phi_{0}\left(
m\right)  -tA\left(  m\right)  \right\}  .
\end{equation}
Here $\Phi_{0}\left(  m\right)  $ is the Legendre-Fenchel transform of
$\Psi\left(  t,x\right)  ,$ evaluated at the initial time $t=0$. \ The
geometric interpretation is that $\Phi\left(  t,m\right)  $ forms the convex
hull of $\Phi^{0}\left(  t,m\right)  $ on the interval $m\in\left[
0,1\right]  $.

\begin{figure}
[h]
\begin{center}
\includegraphics[width=\linewidth]{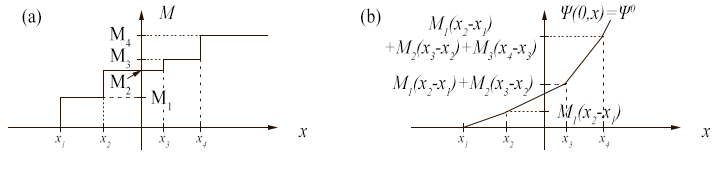}
\caption{(a)\ Cumulative distribution function of mass as a function of
position; (b) Construction of $\Psi\left(  0,x\right)  $.}
\label{Potentialfigure}
\end{center}
\end{figure}

To illustrate this construction, consider a discrete example with mass density
function $M$ along and corresponding function $\Psi$ as graphed in Figure
\ref{Potentialfigure}.

\bigskip%

This corresponds to our example of discrete point masses $m_{i}$ at positions
$x_{i}$. \ Clearly the function $\Psi$ is continuous but not $C^{2}$ or even
$C^{1}$. \ To show that $\Psi$ is a viscosity solution, we need to show that,
given $\phi\in C^{2}\left(  \mathbb{R\times}\left[  0,\infty\right]  \right)
$ such that $\phi\left(  x_{0},t\right)  =\Psi\left(  x_{0},t\right)  $, the
following hold:

(i) $\Psi$ is continuous (this is trivial).

(ii) $\phi\geq\Psi$ in a neighborhood of $x_{0}$ implies $\partial_{t}%
\phi+A\left(  \partial_{x}\phi\right)  \leq0$.

(iii) $\phi\leq\Psi$ in a neighborhood of $x_{0}$ implies $\partial_{t}%
\phi+A\left(  \partial_{x}\phi\right)  \geq0$.

Notationally, we will use a prime to denote the $x$ derivative of $\phi$, i.e.
$\phi^{\prime}\left(  x_{0},0\right)  =\frac{d}{dx}\phi\left(  x_{0},0\right)
,$ for the remainder of this section.

To satisfy (ii), we let $\phi$ be a $C^{2}$ function which lies above $\Psi$.
\ However, we argue that choosing such a function is impossible. \ Observe
that for $\phi\geq\Psi$ to be satisfied in a neighborhood of the point $x_{2}%
$, we need $\phi^{\prime}\leq M_{1}$ in a one-sided neighborhood of $x_{2}$
("left" of $x_{2}$). \ For if not, then $\phi^{\prime}>M_{1}$ in $\left(
x_{2}-\varepsilon,x_{2}\right)  $, some $\varepsilon>0$, and clearly%
\begin{equation}
\phi\left(  x_{2}-\varepsilon,0\right)  <\phi\left(  x_{2},0\right)
-\varepsilon M_{1}=\Psi\left(  x_{2}-\varepsilon,0\right)  ,
\end{equation}
violating our assumption. \ Thus, $\phi^{\prime}\leq M_{1}$ in a neighborhood
$\left(  x_{2}-\varepsilon\,,x_{2}\right)  $.

Similarly, one can show that $\phi^{\prime}\geq M_{2}$ in a neighborhood
$\left(  x_{2},x_{2}+\tilde{\varepsilon}\right)  ,$ $\tilde{\varepsilon}>0$.
\ But since $M_{1}<M_{2}$, this implies $\phi^{\prime}$ is discontinuous at
$x_{2}$, so that $\phi^{\prime\prime}\left(  x_{2},0\right)  $ does not exist.
\ Hence it is impossible to find a $C^{2}$ function $\phi$ satisfying
$\phi\geq\Psi$ in a neighborhood of $x_{2}$.

In the case of (iii), it is indeed plausible that one can find such a $\phi$.
\ Then one wishes to verify that%
\begin{equation}
\left[  \partial_{t}\phi+A\left(  \partial_{x}\phi\right)  \right]  _{\left(
x_{0},0\right)  }\geq0\text{.}%
\end{equation}
Let $\phi$ be given such that $\phi\in C^{2}\left(  \mathbb{R\times}%
[0,\infty)\right)  $, $\phi\left(  x_{0},0\right)  =\Psi\left(  x_{0}%
,0\right)  $, and $\phi\leq\Psi$ in a neighborhood of $x_{0}$. By considering
the cases (i) $x_{0}\not \in \left\{  x_{i}\right\}  _{i=1}^{4}$ and (ii)
$x_{0}$ taking one of the values $\left\{  x_{i}\right\}  _{i=1}^{4}$, one can
show condition (iii) holds. The algebraic details are not particularly
enlightening and are hence omitted.

Similar results relating to ballistic aggregation of particles are studied in
\cite{V}.

\subsection{Returning to the discrete example}

Let us return to our discrete example with four point masses. \ Recall in this
setting that $M_{i}:=\sum_{j=1}^{i}m_{i}$. \ We sketch $M\left(  0,x\right)  $
in Fig. 2(a) and $\Psi^{0}\left(  x\right)  =\Psi\left(  0,\cdot\right)
=\int_{-\infty}^{x}M\left(  t,y\right)  dy$ in Fig. 2(b). \ Clearly, $\Psi
^{0}\left(  m\right)  $ is given by%
\begin{equation}
\Psi^{0}\left(  x\right)  =\left\{
\begin{array}
[c]{c}%
0\\
M_{1}\left(  x-x_{1}\right) \\
M_{1}\left(  x_{2}-x_{1}\right)  +M_{2}\left(  x-x_{2}\right) \\
M_{1}\left(  x_{2}-x_{1}\right)  +M_{2}\left(  x_{3}-x_{2}\right)
+M_{3}\left(  x-x_{3}\right) \\
M_{1}\left(  x_{2}-x_{1}\right)  +M_{2}\left(  x_{3}-x_{2}\right)
+M_{3}\left(  x_{4}-x_{3}\right)  +M_{4}x
\end{array}
\right.
\begin{array}
[c]{c}%
x<x_{1}\\
x_{1}<x<x_{2}\\
x_{2}<x<x_{3}\\
x_{3}<x<x_{4}\\
x_{4}<x
\end{array}
\label{psipiecewise}%
\end{equation}
It is convenient to express $\Psi^{0}\left(  x\right)  $ in terms of positive
parts of functions rather then the piecewise construction in
(\ref{psipiecewise}), i.e. we write
\begin{equation}
\Psi^{0}\left(  x\right)  =\sum_{i=1}^{4}M_{1}\left(  x-x_{i}\right)
_{+}=m_{1}\left(  x-x_{1}\right)  _{+}+m_{2}\left(  x-x_{2}\right)  _{+}%
+m_{3}\left(  x-x_{3}\right)  _{+}+m_{4}\left(  x-x_{4}\right)  _{+}%
\end{equation}
This is easy to see intuitively (for example, the $M_{1}\left(  x-x_{1}%
\right)  $ term becomes $m_{1}\left(  x_{2}-x_{1}\right)  $ for $x>x_{2}$, and
the $M_{2}\left(  x-x_{2}\right)  $ contributes the $m_{1}\left(
x-x_{2}\right)  $ part). \ Now, we compute $\Phi^{0}\left(  m\right)
=\Psi^{\ast}\left(  0,m\right)  =\sup_{x\in\mathbb{R}}\left\{  xm-\Psi\left(
0,x\right)  \right\}  .$ \ We differentiate the expression and find:%
\begin{equation}
Q_{x}:=\left(  xm-\Psi\left(  0,x\right)  \right)  _{x}=m-m_{1}H\left(
x-x_{1}\right)  -m_{2}H\left(  x-x_{2}\right)  -m_{3}H\left(  x-x_{3}\right)
-m_{4}H\left(  x-x_{4}\right)  .
\end{equation}
Note: (i) For $m\not =m_{i}$, it has no critical points, (ii) for any\ value
of $m$, $Q_{x}$ is a decreasing function. \ For example, for $m<m_{1}$,
clearly the maximum is attained at $x=x_{1}$, and hence%
\begin{align}
\Phi^{0}\left(  m\right)   & =x_{1}m-\Psi\left(  0,x_{1}\right)  =x_{1}%
m-m_{1}\left(  x_{1}-x_{1}\right)  _{+}-m_{2}\left(  x_{1}-x_{2}\right)
_{+}-m_{3}\left(  x_{1}-x_{3}\right)  _{+}-m_{4}\left(  x_{1}-x_{4}\right)
_{+}\nonumber\\
& =x_{1}m
\end{align}
After performing similar computations for other cases and combining the
results, one is lead to%
\begin{equation}
\Phi^{0}\left(  m\right)  =x_{1}m+\sum_{i=1}^{n-1}\left(  m-\sum_{j=1}%
^{i}m_{j}\right)  _{+}\left(  x_{i+1}-x_{i}\right)  .
\end{equation}
To justify (29) in the paper, we consider the discrete case as in Fig. 2(b)
and observe%
\begin{align}
A\left(  m\right)   & =\left\{
\begin{array}
[c]{c}%
v_{1}m\\
v_{1}m_{1}+v_{2}\left(  m-m_{1}\right) \\
v_{1}m_{1}+v_{2}m_{2}+v_{3}\left(  m-M_{2}\right) \\
v_{1}m_{1}+v_{2}m_{2}+v_{3}m_{3}+v_{4}\left(  m-M_{3}\right)
\end{array}
\right.
\begin{array}
[c]{c}%
m<M_{1}\\
M_{1}<m<M_{2}\\
M_{2}<m<M_{3}\\
M_{3}<m<M_{4}%
\end{array}
\nonumber\\
& =v_{1}m+\sum_{i=1}^{3}\left(  m-\sum_{j=1}^{i}m_{i}\right)  \left(
v_{i+1}-v_{i}\right)
\end{align}
using the same arguments as before for rearranging the function. \ Both
arguments carry over easily using induction for the algebra.

\begin{figure}
[h]
\begin{center}
\includegraphics[width=\linewidth]{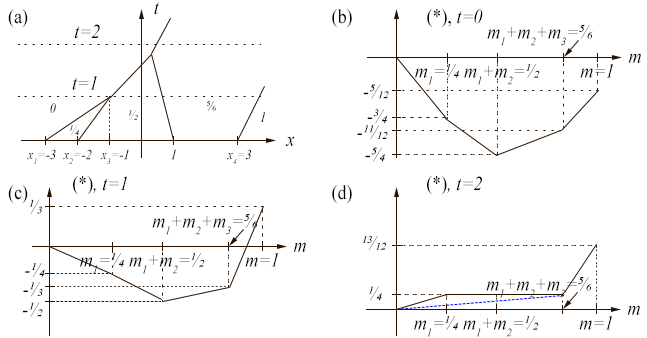}
\caption{(a) Representation of mass in cumulative distribution form up to a
point $x$ in the $xt$ plane; (b)-(d)\ Plot of the expression $\Phi^{0}\left(
m\right)  +tA\left(  m\right)  $ for times $t=0,1,2$ respectively in solid
lines; following the dashed lines forms the convex hull, yielding the Legendre
transform $\Phi_{n}\left(  t,m\right)  .$ Note that for (b) and (c), the
expression and its convex hull are identical, and in (d) there is a
distinction, with the convex hull indicated by the dashed blue line.}
\label{Masscdffigure}
\end{center}
\end{figure}

Then the expression%
\begin{equation}
\Phi_{n}\left(  t,m\right)  =x_{1}\left(  t\right)  m+\sum_{i=1}^{n-1}\left(
m-\sum_{j=1}^{i}m_{j}\right)  _{+}\left(  x_{i+1}\left(  t\right)
-x_{i}\left(  t\right)  \right)
\end{equation}
is just given by the convex hull of%
\begin{equation}
\Phi^{0}\left(  m\right)  +tA\left(  m\right)  =\left(  x_{1}+tv_{1}\right)
m+\sum_{i=1}^{n-1}\left(  m-\sum_{i=1}^{j}m_{i}\right)  _{+}\left(
x_{i+1}+tv_{i+1}-x_{i}-tv_{i}\right) \label{BGpotential}%
\end{equation}
To illustrate this, we return to our example. \ Let%
\begin{align}
m_{1}  & =\frac{1}{4},\text{ \ }m_{2}=\frac{1}{4},\text{ \ }m_{3}=\frac{1}%
{3},\text{ \ }m_{4}=\frac{1}{6}\nonumber\\
v_{1}  & =2,\text{ \ }v_{2}=1,\text{ \ }v_{3}=-\frac{1}{2},\text{ \ }%
v_{4}=1\nonumber\\
x_{1}  & =-3,\text{ \ }x_{2}=-2,\text{ \ }x_{3}=1,\text{ \ }x_{4}=3
\end{align}
Note $\sum_{i}m_{i}=1$. \ One can readily compute the quantity
(\ref{BGpotential}) by breaking it into cases for the values of $m$, and has%
\begin{equation}
LHS\text{ of (\ref{BGpotential})}=\left\{
\begin{array}
[c]{c}%
\left(  -3+2t\right)  m\\
\left(  -2+t\right)  m-\frac{1}{4}\left(  1-t\right) \\
\left(  1-\frac{t}{2}\right)  m-\frac{1}{4}\left(  7-4t\right) \\
\left(  3+t\right)  m-\frac{1}{12}\left(  41+3t\right)
\end{array}
\right.
\begin{array}
[c]{c}%
m<\frac{1}{4}\\
\frac{1}{4}<m<\frac{1}{2}\\
\frac{1}{2}<m<\frac{5}{6}\\
\frac{5}{6}<m<1
\end{array}
\end{equation}
Now we want to calculate (*) for different values of $t$ and illustrate what
the convex hull looks like. \ By simple geometry, one observes that particles
$1$ and $2$ will collide and stick at time $t=1$, and subsequently that will
collide with $3$ at some time $1<t<2$ (see Figure \ref{Masscdffigure}(a),
where the double-dashed line represents $t=2$). Note that although this is
defined piecewise, it is continuous since the boundary values match at
$\frac{1}{4},$ $\frac{1}{2},$ and $\frac{5}{6}$. We consider (*) for three
different values of $t$.

For the initial state, $t=0$, one has%
\begin{equation}
LHS\text{ of (\ref{BGpotential})}|_{t=0}=\left\{
\begin{array}
[c]{c}%
-3m\\
-2m-\frac{1}{4}\\
m-\frac{7}{4}\\
3m-\frac{41}{12}%
\end{array}
\right.
\begin{array}
[c]{c}%
m<\frac{1}{4}\\
\frac{1}{4}<m<\frac{1}{2}\\
\frac{1}{2}<m<\frac{5}{6}\\
\frac{5}{6}<m<1
\end{array}
\end{equation}

Similarly, for $t=1,2$, we have, respectively%
\begin{equation}
LHS\text{ of (\ref{BGpotential})}||_{t=1}=\left\{
\begin{array}
[c]{c}%
-m\\
-m\\
\frac{m}{2}-\frac{3}{4}\\
4m-\frac{11}{3}%
\end{array}
\right.
\begin{array}
[c]{c}%
m<\frac{1}{4}\\
\frac{1}{4}<m<\frac{1}{2}\\
\frac{1}{2}<m<\frac{5}{6}\\
\frac{5}{6}<m<1
\end{array}
\end{equation}%
\begin{equation}
LHS\text{ of (\ref{BGpotential})}||_{t=2}=\left\{
\begin{array}
[c]{c}%
m\\
\frac{1}{4}\\
\frac{1}{4}\\
5m-\frac{47}{12}%
\end{array}
\right.
\begin{array}
[c]{c}%
m<\frac{1}{4}\\
\frac{1}{4}<m<\frac{1}{2}\\
\frac{1}{2}<m<\frac{5}{6}\\
\frac{5}{6}<m<1
\end{array}
\end{equation}
which is once again continuous.%

These are graphed in Figure \ref{Masscdffigure}(b), (c), and (d),
respectively. \ Recall that we wanted to consider the convex hull of the
expression above. \ In the first two cases, the resulting function is already
convex. \ In the third case (after shocks have occurred), the function is not
convex and the convex hull is formed by taking the piecewise linear function
given by the double-dashed line for $0<m<\frac{5}{6}$ and the piece
$5m-\frac{47}{12}$ for $\frac{5}{6}<m<1$.

\section{Analysis of Burgers' Equation Using Flow Maps}

\subsection{Definition of weak solution for conservation laws}

Another approach to this problem is to use a Generalised Variational Principle
(GVP) as in E, Rykov, and Sinai \cite{ERS}. \ The first task is to formulate
the equations in a weak form. \ Specifically, consider the system of
conservation laws%
\begin{align}
\rho_{t}+\left(  \rho u\right)  _{x}  & =0\nonumber\\
\left(  \rho u\right)  _{t}+\left(  \rho u^{2}\right)  _{x}  & =0.\label{fmic}%
\end{align}
We want to define a weak solution by having (\ref{fmic}) hold when we multiply
by a test function and integrate. \ Since $\rho$ is a purely singular measure
for any time $t>0$, it is incorrect to write%
\begin{align}
\int\rho\phi_{t}+\rho u\phi_{x}  & =0+\text{ boundary terms}\nonumber\\
\int\rho u\phi_{t}+\rho u^{2}\phi_{x}  & =0+\text{ boundary terms}%
\end{align}
Instead, we take a family $\left(  P_{t},I_{t}\right)  $ of Borel measures
that are weakly continuous with respect to $t$ such that $I_{t}$ is absolutely
continuous with respect to $P_{t}$ for each fixed $t$. \ Then we can define
the Radon-Nikodym derivative as%
\begin{equation}
u\left(  x,t\right)  =\frac{dI_{t}}{dP_{t}}\left(  x\right)  ,
\end{equation}
i.e.%
\begin{equation}
\int_{A}u\left(  x,t\right)  dP_{t}\left(  x\right)  =\int_{A}dI_{t}\left(
x\right)  .
\end{equation}
for any measurable set $A$ in the appropriate set of functions. \ More
specifically,%
\begin{equation}
\int u\left(  x,t\right)  f\left(  x\right)  dP_{t}\left(  x\right)  =\int
f\left(  x\right)  dI_{t}\left(  x\right)
\end{equation}
for any measurable function $f$. We can now integrate and call $\left(
P_{t},I_{t},u\right)  _{t\geq0}$ a weak solution of (\ref{fmic}) if it
satisfies the resulting equality. \ More precisely, for any $f,g\in C_{0}%
^{1}\left(  \mathbb{R}\right)  $ and $0<t_{1}<t_{2}$, we need%
\begin{align}
& \int\int_{t_{1}}^{t_{2}}\left\{  \rho_{t}\left(  \eta\right)  f\left(
\eta\right)  +\left(  \rho u\right)  _{x}f\right\}  d\tau d\eta\nonumber\\
& =\int f\left(  \eta\right)  dP_{t_{2}}\left(  \eta\right)  -\int f\left(
\eta\right)  dP_{t_{1}}\left(  \eta\right)  -\int_{t_{1}}^{t_{2}}d\tau\int
f^{\prime}\left(  \eta\right)  u\left(  \eta,\tau\right)  dP_{\tau}\left(
\eta\right) \nonumber\\
& =\int f\left(  \eta\right)  dP_{t_{2}}\left(  \eta\right)  -\int f\left(
\eta\right)  dP_{t_{1}}\left(  \eta\right)  -\int_{t_{1}}^{t_{2}}d\tau\int
f^{\prime}\left(  n\right)  dI_{\tau}\left(  \eta\right)  =0,
\end{align}
where the boundary terms from the integration by parts in the second term drop
out since $f$ has compact support. \ Note that there is no integration by
parts in evaluating the first term ($f$ has no dependence on $t$ so we just
integrate $\rho_{t}$). \ To write the second equation of (\ref{fmic}) in weak
form, we proceed similarly:%
\begin{align}
& \int\int_{t_{1}}^{t_{2}}\left\{  \left(  \rho\left(  \eta\right)  u\right)
_{t}g\left(  \eta\right)  +\left(  \rho\left(  \eta\right)  u^{2}\right)
_{x}g\left(  \eta\right)  \right\}  d\tau d\eta\nonumber\\
& =\int g\left(  \eta\right)  u\left(  \eta,t_{2}\right)  dP_{t_{2}}\left(
\eta\right)  -\int g\left(  \eta\right)  u\left(  \eta,t_{1}\right)
dP_{t_{1}}\left(  \eta\right) \nonumber\\
& -\int_{t_{1}}^{t_{2}}d\tau\int g^{\prime}\left(  \eta\right)  u\left(
\eta,\tau\right)  u\left(  \eta,\tau\right)  dP_{\tau}\left(  \eta\right)
\nonumber\\
& =\int g\left(  \eta\right)  dI_{t_{2}}\left(  \eta\right)  -\int g\left(
\eta\right)  dI_{t_{1}}\left(  \eta\right)  -\int_{t_{1}}^{t_{2}}d\tau\int
g^{\prime}\left(  \eta\right)  u\left(  \eta,\tau\right)  dI_{\tau}\left(
\eta\right)  =0
\end{align}
Rewriting, our definition of a weak solution under the above assumptions
becomes:%
\begin{align}
\int f\left(  \eta\right)  dP_{t_{2}}\left(  \eta\right)  -\int f\left(
\eta\right)  dP_{t_{1}}\left(  \eta\right)   & =\int_{t_{1}}^{t_{2}}d\tau\int
f^{\prime}\left(  n\right)  dI_{\tau}\left(  \eta\right) \nonumber\\
\int g\left(  \eta\right)  dI_{t_{2}}\left(  \eta\right)  -\int g\left(
\eta\right)  dI_{t_{1}}\left(  \eta\right)   & =\int_{t_{1}}^{t_{2}}d\tau\int
g^{\prime}\left(  \eta\right)  u\left(  \eta,\tau\right)  dI_{\tau}\left(
\eta\right)
\end{align}
which matches Definition 1 in \cite{ERS}.

\subsection{Discrete example}

We again consider our discrete example with four particles, with initial
conditions specified by%
\begin{align}
m_{1}^{0}  & =\frac{1}{4},\text{ \ }m_{2}^{0}=\frac{1}{4},\text{ \ }m_{3}%
^{0}=\frac{1}{3},\text{ \ }m_{4}^{0}=\frac{1}{6},\nonumber\\
v_{1}^{0}  & =2,\text{ \ }v_{2}^{0}=1,\text{ \ }v_{3}^{0}=-\frac{1}{2},\text{
\ }v_{4}^{0}=1,\nonumber\\
x_{1}^{0}  & =-3,\text{ \ }x_{2}^{0}=-2,\text{ \ }x_{3}^{0}=1,\text{ \ }%
x_{4}^{0}=3,
\end{align}
as before. \ We first want to compute the flow map under the following initial
velocity field:%
\begin{equation}
u_{0}\left(  x\right)  =\left\{
\begin{array}
[c]{c}%
0\\
2\\
1\\
-\frac{1}{2}\\
1
\end{array}
\right.
\begin{array}
[c]{c}%
x<-3\\
-3\leq x<-2\\
-2\leq x<1\\
1\leq x<3\\
3\leq x
\end{array}
\end{equation}
An equivalent method would be to consider a velocity field with $\delta
$-functions. at each point mass, but this is omitted for brevity. \ For time
$t$ before the first collision at $\left(  x_{12},t_{1}^{\ast}\right)
=\left(  -1,1\right)  $, the flow map is as follows:%
\begin{equation}
\varphi_{t}\left(  x\right)  =\left\{
\begin{array}
[c]{c}%
x\\
-3+2t\\
x+2t\\
-2+t\\
x+t\\
x-\frac{t}{2}\\
3+t\\
x+t
\end{array}
\right.
\begin{array}
[c]{c}%
x<-3\\
x=-3\\
-3<x<-2\\
x=-2\\
-2<x<1\\
1\leq x<3\\
x=3\\
3<x
\end{array}
\end{equation}
Therefore, the inverse flow map is given by%
\begin{equation}
\varphi_{t}^{-1}\left(  x\right)  =\left\{
\begin{array}
[c]{c}%
x\\
\left\{  \emptyset\right\} \\
\left[  -3,-3+2t\right] \\
x\\
\left\{  \emptyset\right\} \\
\left[  -2,-2+t\right] \\
x\\
\left[  1-\frac{t}{2},1\right] \\
\left\{  \emptyset\right\} \\
x\\
\left\{  \emptyset\right\} \\
\left[  3,3+t\right] \\
x
\end{array}
\right.
\begin{array}
[c]{c}%
x<-3\\
-3\leq x<-3+2t\\
x=-3+2t\\
-3+2t<x<-2\\
-2\leq x<-2+t\\
x=-2+t\\
-2+t<x<1\\
x=1-\frac{t}{2}\\
1-\frac{t}{2}<x\leq1\\
1<x<3\\
3\leq x<3+t\\
x=3+t\\
3+t<x
\end{array}
\end{equation}
The elements of the partition $\xi_{t}$ are then given by real numbers in the
intervals $\left(  -\infty,-2+t\right)  ,$ $\left(  -2+2t,1-\frac{t}%
{2}\right)  ,$ and $\left(  1+t,\infty\right)  $ along with the intervals
$\left[  -2+2t,1-\frac{t}{2}\right]  $ and $\left[  1-\frac{t}{2},1+t\right]
.$ \ As we increase time, these intervals will grow and eventually merge when
we have collisions between 1 and 2, and then the resulting particle with 3.

\begin{figure}
[h]
\begin{center}
\includegraphics[width=\linewidth]{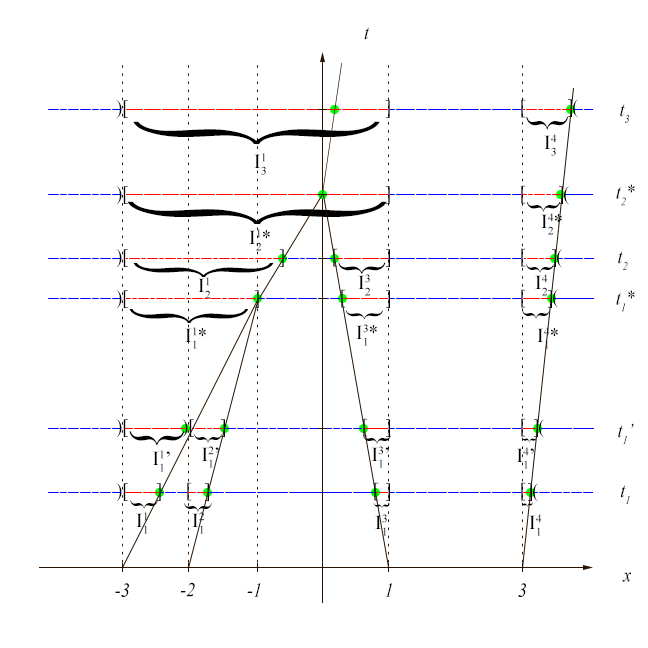}
\caption{Evolution of the discrete example and mapping back using the flow
map. Highlighted in blue (long-short dash lines) are intervals unchanged under
the flow map. In red (long dashed line) are intervals for which the flow map
inverse is undefined. The points in green correspond to single points for
which an entire interval is mapped back onto, which occurs in notably many
cases. For example, $\varphi_{t_{2}^{\ast}}^{-1}\left(  I_{2}^{1\ast}\right)
=\left\{  0\right\}  $ and $\varphi_{t_{1}^{\ast}}^{-1}\left(  \left\{
-2\right\}  \right)  =\left\{  \emptyset\right\}  $.}
\label{Flowmapfigure}
\end{center}
\end{figure}

From here we can identify the element $C_{t}\left(  y\right)  $ corresponding
to the element of the partition $\xi_{t}$ containing $y$, and reconstruct the
solution using equation (1.10) in the paper:%
\begin{equation}
\varphi_{t}\left(  y\right)  =\frac{\int_{C_{t}\left(  y\right)  }\left(
\eta+tu_{0}\left(  \eta\right)  \right)  dP_{0}\left(  \eta\right)  }%
{\int_{C_{t}\left(  y\right)  }dP_{0}\left(  \eta\right)  },\text{
\ \ }u\left(  x,t\right)  =\frac{\int_{D_{t}\left(  x\right)  }u_{0}\left(
\eta\right)  dP_{0}\left(  \eta\right)  }{\int_{D_{t}\left(  x\right)  }%
dP_{0}\left(  \eta\right)  }.
\end{equation}
Construction of the inverse flow map is shown in Figure \ref{Flowmapfigure}.
Note that our calculations are formal, but can be made rigourous through
application of the lemmas and Theorem in \cite{ERS}. In a similar vein, one
can also apply front tracking methods as described in \cite{HO}.

\subsection{GVP and the continuous case}

We give a brief idea of generalisation to the continuous case. \ To do this,
we will need several assumptions. \ We provide both the technical definition
and an explanation of the physical meaning of each. \ We first let
$P_{0},I_{0}\in M$, the space of Radon measures on $\mathbb{R}^{1}$,
$P_{0}\geq0$. \ A Radon measure is defined as a measure that is inner regular
(for all Borel sets $B$, $m\left(  B\right)  =\sup\left\{  m\left(  K\right)
|K\subset B\text{ compact}\right\}  $) measure defined on the $\sigma$-algebra
of a Hausdorff topological space $X$ that is locally finite (for every point
of $X$, there exists a neighborhood $U$ such that $m\left(  U\right)  <\infty$).

\textbf{(Assumption 1).} For any compact $\Lambda\subset\mathbb{R}^{1},$ one
has $P_{0}\left(  \Lambda\right)  <\infty$ and $P_{0}$ is either discrete or
absolutely continuous with respect to Lebesgue measure. \ If $P_{0}~$is
absolutely continuous, then we assume $\rho_{0}\left(  x\right)  >0$ for all
points $x$ in the support of $\rho_{0}$. \ If $Supp\left(  \rho_{0}\right)  $
is unbounded, then we assume%
\begin{equation}
\int_{0}^{x}sdP_{0}\left(  s\right)  \rightarrow+\infty\text{ as }\left\vert
x\right\vert \rightarrow\infty\text{.}%
\end{equation}
The first statement corresponds to the physical requirement that we do not
have an infinite amount of mass on any finite interval (intervals are
precompact in $\mathbb{R}$), and that the distribution of masses must be
either (i) fully discrete, with point masses $\left\{  m_{i}\right\}  $ at
locations $\left\{  x_{i}\right\}  $, or (ii) there are not any such point
masses anywhere on $\mathbb{R}$. \ In this second case, we then assume it has
a density, with either a finite cutoff, or density out to infinity that does
not fall off too sharply (less sharply than $\frac{1}{s^{2+\varepsilon}}$, for example).

\textbf{(Assumption 2).} \ The initial distribution of momentum $I_{0}$ is
absolutely continuous with respect to $P_{0}$. \ We can define a Radon-Nikodym
derivative $u\left(  \cdot,0\right)  =\frac{dI_{0}}{dP_{0}}$ and this is the
initial velocity. \ When $P_{0}$ is absolutely continuous, we assume further
that $u\left(  \cdot,0\right)  $ is continuous as well. \ This corresponds
simply to the requirement that momentum is zero on intervals where there is no mass.

\textbf{(Assumption 3).} \ For all $z>0$, we have%
\begin{equation}
\sup_{\left\vert x\right\vert \leq z}\left\vert u_{0}\left(  x\right)
\right\vert \leq b_{0}\left(  z\right)  \text{ where }\lim_{\left\vert
z\right\vert \rightarrow\infty}\frac{1}{z}b_{0}\left(  z\right)  =0.
\end{equation}
Physically this means that the initial velocities of the particles can not
increase as $O\left(  \left\vert z\right\vert \right)  $ or more as position
goes to infinity.

Under these assumptions, we have a way of constructing the partition $\xi_{t}$
using the initial data, using a Generalised Variational Principle. \ We have
that $y\in\mathbb{R}^{1}$ is a left endpoint of an element of $\xi_{t}$ if and
only for every $y^{-},y^{+}\in\mathbb{R}$ such that $y^{-}<y<y^{+}$, we have
the following:%
\begin{equation}
\frac{\int_{\lbrack y^{-},y)}\left(  \eta+tu\left(  \eta;0\right)  \right)
dP_{0}}{\int_{[y^{-},y)}dP_{0}\left(  \eta\right)  }<\frac{\int_{[y,y^{+}%
]}\left(  \eta+tu\left(  \eta;0\right)  \right)  dP_{0}\left(  \eta\right)
}{\int_{[y,y^{+}]}dP_{0}\left(  \eta\right)  }.
\end{equation}

\section{Entropy Solution and Variational Approach}

We now consider the work of Huang and Wang \cite{HW} on the system of
one-dimensional pressureless gas equations given by%
\begin{align}
\rho_{t}+\left(  \rho u\right)  _{x}  & =0\label{HWstart}\\
\left(  \rho u\right)  _{t}+\left(  \rho u^{2}\right)  _{x}  & =0.\nonumber
\end{align}
A basis of their approach entails generalising characteristics when the flow
map breaks down (is no longer one-to-one). \ In particular, they consider the
generalised potential given by%
\begin{equation}
F\left(  y;x,t\right)  =\int_{0+0}^{y-0}\left(  tu_{0}\left(  \eta\right)
+\eta-x\right)  dm_{0}\left(  \eta\right) \label{HWgp}%
\end{equation}
and first note that if $\rho,u$ are bounded and measurable functions, then
$m\left(  x,t\right)  =%
{\displaystyle\oint\nolimits_{\left(  0,0\right)  }^{\left(  x,t\right)  }}
\rho dx-\rho udt$ is independent of path since $\rho_{t}+\left(  \rho
u\right)  _{x}$ is conserved. \ Indeed, $\rho_{t}=-\left(  \rho u\right)
_{x}=\left(  -\rho u\right)  _{x}$. \ Further, we have $m_{x}=\rho,m_{t}=-\rho
u$. \ One then readily verifies that (\ref{HWstart}) is then equivalent to%
\begin{align}
m_{t}+um_{x}  & =0\nonumber\\
\left(  m_{x}u\right)  _{t}+\left(  m_{x}u^{2}\right)  _{x}  &
=0\label{HWsystem}%
\end{align}
A weak solution to the system (\ref{HWsystem}) is defined as follows.

\begin{definition}
Let $m\left(  x,t\right)  $ be of bounded variation locally in $x$, and
$u\left(  x,t\right)  $ be bounded and $m_{x}$-measurable. \ Assume
$m_{x},um_{x}$ are weakly continuous in $t$. \ We call $\left(  \rho,u\right)
=\left(  m_{x},u\right)  $ a \textit{weak solution} of (\ref{HWsystem}) given
that%
\begin{align}
\int\int\varphi_{t}mdxdt-\int\int\varphi udmdt  & =0\nonumber\\
\int\int\psi_{t}u+\psi_{x}u^{2}dmdt  & =0
\end{align}
is satisfied for all $\varphi,\psi\in C_{0}^{\infty}\left(  \mathbb{R}_{+}%
^{2}\right)  ,$ where the integrals are the Lebesgue-Stieltjes integrals.
\end{definition}

We understand the initial value in the sense that as we take $\eta$ to the
lower limit, the measures $\rho,\rho u$ converge weakly (\cite{BZ}, p. 57),
where $u_{0}$ is bounded and is measurable with respect to $\rho_{0}$. \ We
also have the following entropy condition: we call $\left(  \rho,u\right)  $
an entropy condition if%
\begin{equation}
\frac{u\left(  x_{2},t\right)  -u\left(  x_{1},t\right)  }{x_{2}-x_{1}}%
\leq\frac{1}{t}%
\end{equation}
holds for any $x_{1}<x_{2}$, a.e. in $t>0$, and $\rho u^{2}$ converges weakly
to $\rho_{0}u_{0}^{2}$ as $t\rightarrow0$.

Their main result is the following:

\begin{theorem}
(Existence). \ Let $\rho_{0}\geq0\in M_{loc}\left(  \mathbb{R}\right)  $, the
space of Radon measures defined on $\mathbb{R}$ (or $m_{0}\left(  x\right)
=\rho([0,x))$ increasing) and let $u_{0}$ be bounded and measurable with
respect to $\rho_{0}$, then the system (\ref{HWsystem}) admits at least one
entropy solution.
\end{theorem}

In proving the result, a number of technical lemmas are required. By using the
generalised potential, one uses them to construct an entropy solution for an
increasing function $m_{0}\left(  x\right)  =\rho_{0}\left(  [0,x)\right)  $.
\ The trivial case $\rho_{0}=0$ is excluded. \ We state several of the lemmas
without proof; more details can be found in \cite{HW}.

\bigskip

\textbf{Lemma 1. \ }For any point $\left(  x,t\right)  $, the function
$F\left(  y;x,t\right)  $ considered as a function of $y$ has a finite lower bound.

Now we define
\begin{align}
v\left(  x,t\right)   & =\min_{y}F\left(  y;x,t\right)  ,\nonumber\\
S\left(  x,t\right)   & =\left\{  y|\text{there exists }y_{n}\rightarrow
y\text{ s.t. }F\left(  y_{n};x,t\right)  \rightarrow v\left(  x,t\right)
\right\}  ,
\end{align}
i.e. $S\left(  x,t\right)  $ as the set of points $y$ for which we can find a
sequence $y_{n}$ approaching this limit. \ Using the fact that $F$ is left
continuous in $y,$ for every $y_{0}\in S\left(  x,t\right)  $, one has%
\begin{equation}
v\left(  x,t\right)  =\left\{
\begin{array}
[c]{c}%
F\left(  y_{0};x,t\right) \\
F\left(  y_{0}+0;x,t\right)
\end{array}
\right.
\begin{array}
[c]{c}%
F\left(  y;x,t\right)  \text{ achieves its minimum at }y_{0}\\
\text{otherwise}%
\end{array}
,
\end{equation}
leading to the following lemma.

\bigskip

\textbf{Lemma 2. \ }Assume $y_{0}\in S\left(  x,t\right)  $, $\left[
m_{0}\left(  y_{0}\right)  \right]  =m_{0}\left(  y_{0}+0\right)
-m_{0}\left(  y_{0}-0\right)  >0$. \ Then%
\begin{equation}
v\left(  x,t\right)  =\min_{y}F\left(  y;x,t\right)  =\left\{
\begin{array}
[c]{c}%
F\left(  y_{0};x,t\right) \\
F\left(  y_{0}+0;x,t\right)
\end{array}
\right.
\begin{array}
[c]{c}%
if\text{ }x\leq y_{0}+tu_{0}\left(  y_{0}\right) \\
if\text{ }x>y_{0}+tu_{0}\left(  y_{0}\right)
\end{array}
\end{equation}
In other words, in the first case, $F$ achieves its minimum, and in the second
it does not.

\bigskip

\textbf{Lemma 3. \ }Let $\left(  x_{n},t_{n}\right)  $ and $y_{n}\in S\left(
x_{n},t_{n}\right)  $ converge to $\left(  x,t\right)  $ and $y_{0}$,
respectively. \ Then $y_{0}\in S\left(  x,t\right)  $. From Lemma 2.1, we have
that $\inf\left\{  y|y\in spt\left\{  \rho_{0}\right\}  \right\}  $ is finite
if $y_{m}\left(  x,t\right)  =-\infty$. \ Similarly, $\sup\left\{  y|y\in
spt\left\{  \rho_{0}\right\}  \right\}  $ is finite if $y^{m}\left(
x,t\right)  =+\infty.$ For each point $\left(  x_{0},t_{0}\right)  $,
introduce left and right backward generalised characteristics $L_{1},L_{2}:$%
\begin{align}
L_{1}  & :x=x_{0}+\frac{x_{0}-y_{\ast}\left(  x_{0},t_{0}\right)  }{t_{0}%
}\left(  t-t_{0}\right) \nonumber\\
L_{2}  & :x=x_{0}+\frac{x_{0}-y^{\ast}\left(  x_{0},t_{0}\right)  }{t_{0}%
}\left(  t-t_{0}\right)
\end{align}
and claim that there is only one minimum point of $F\left(  y;x,t\right)  $
for each $\left(  x,t\right)  $ along backward lines $L_{1},L_{2}$.

\bigskip

\textbf{Lemma 4. \ }For any $y_{0}\in S\left(  x_{0},t_{0}\right)  ,$
$y_{\ast}\left(  x,t\right)  =\allowbreak y^{\ast}\left(  x,t\right)  $ holds
along the lines%
\begin{equation}
L:x=x_{0}+\frac{x_{0}-y_{0}}{t_{0}}\left(  t-t_{0}\right)  \text{.}%
\end{equation}
Furthermore, $y_{\ast}\left(  x,t\right)  =y^{\ast}\left(  x,t\right)  \leq
y_{\ast}\left(  x_{0},t_{0}\right)  $ along the line $L_{1}$, and $y_{\ast
}\left(  x,t\right)  =y^{\ast}\left(  x,t\right)  =y^{\ast}\left(  x_{0}%
,t_{0}\right)  $ along the line $L_{2}$.

\bigskip

\textbf{Lemma 5. \ }$y_{\ast},y^{\ast}$ are increasingly monotonic in $x$.
\ In particular, $y^{\ast}\left(  x_{1},t\right)  \leq y_{\ast}\left(
x_{2},t\right)  $ holds for any $x_{1}<x_{2}$.

\subsection{Constructing the Generalised Potential For a Discrete Example}

\begin{figure}
[ptb]
\begin{center}
\includegraphics[width=\linewidth]{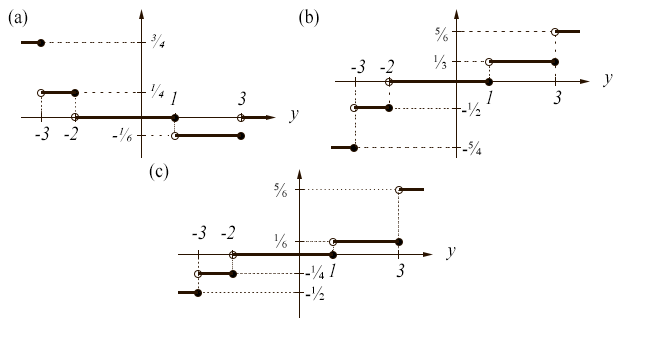}
\caption{Graphs of (a) $\int_{0+0}^{y-0}tu_{0}\left(  \eta\right)
dm_{0}\left(  \eta\right)  $ (with $t=1$), (b) $\int_{0+0}^{y-0}\left(
\eta-x\right)  dm_{0}\left(  \eta\right)  $, (c) $F\left(  y;0,1\right)  $}
\label{HWfigure}
\end{center}
\end{figure}

In conveying the ideas of the construction of the solution and application of
the lemmas, it is useful to return to our discrete example given by four point
masses on the real line, i.e.%
\begin{align}
m_{1}^{0}  & =\frac{1}{4},m_{2}^{0}=\frac{1}{4},m_{3}^{0}=\frac{1}{3}%
,m_{4}^{0}=\frac{1}{6}\nonumber\\
v_{1}^{0}  & =2,v_{2}^{0}=1,v_{3}^{0}=-\frac{1}{2},v_{4}^{0}=1\nonumber\\
x_{1}^{0}  & =-3,x_{2}^{0}=-2,x_{3}^{0}=1,x_{4}^{0}=3
\end{align}
We have $m_{0}\left(  x\right)  =\rho_{0}\left(  [0,x)\right)  .$ \ For
negative $x$ this is interpreted as $m_{0}\left(  x\right)  =-\rho_{0}\left(
(x,0]\right)  .$ \ Hence%
\begin{equation}
m_{0}\left(  x\right)  =\rho_{0}\left(  [0,x)\right)  =\left\{
\begin{array}
[c]{c}%
-m_{1}^{0}-m_{2}^{0}\\
-m_{2}^{0}\\
0\\
m_{3}^{0}\\
m_{3}^{0}+m_{4}^{0}%
\end{array}
\right.
\begin{array}
[c]{c}%
x<x_{1}^{0}\\
x_{1}^{0}\leq x<x_{2}^{0}\\
x_{2}^{0}\leq x\leq x_{3}^{0}\\
x_{3}^{0}<x\leq x_{4}^{0}\\
x_{4}^{0}<x
\end{array}
=\left\{
\begin{array}
[c]{c}%
-\frac{1}{2}\\
-\frac{1}{4}\\
0\\
\frac{1}{3}\\
\frac{1}{6}%
\end{array}
\right.
\begin{array}
[c]{c}%
x<-3\\
-3\leq x<-2\\
-2\leq x\leq1\\
1<x\leq3\\
3<x
\end{array}
\end{equation}
Evaluating the Stieltjes integral for the generalised potential (\ref{HWgp})
for $x=0$ and $t=1,$one has%
\begin{equation}
F\left(  y;0,1\right)  =\left\{
\begin{array}
[c]{c}%
t\left(  v_{1}^{0}m_{1}^{0}+v_{2}^{0}m_{2}^{0}\right)  +\left(  x_{1}%
^{0}-x\right)  m_{1}^{0}+\left(  x_{2}^{0}-x\right)  m_{2}^{0}\\
tv_{2}^{0}m_{2}^{0}+\left(  x_{2}^{0}-x\right)  m_{2}^{0}\\
0\\
tv_{3}^{0}m_{3}^{0}+\left(  x_{3}^{0}-x\right)  m_{3}^{0}\\
t\left(  v_{3}^{0}m_{3}^{0}+v_{4}^{0}m_{4}^{0}\right)  +\left(  x_{3}%
^{0}-x\right)  m_{3}^{0}+\left(  x_{4}^{0}-x\right)  m_{4}^{0}%
\end{array}
\right.
\begin{array}
[c]{c}%
y\leq x_{1}^{0}\\
x_{1}^{0}<y\leq x_{2}^{0}\\
x_{2}^{0}<y\leq x_{3}^{0}\\
x_{3}^{0}<y\leq x_{4}^{0}\\
x_{4}^{0}<y
\end{array}
=\left\{
\begin{array}
[c]{c}%
-\frac{1}{2}\\
-\frac{1}{4}\\
0\\
\frac{1}{6}\\
\frac{5}{6}%
\end{array}
\right.
\begin{array}
[c]{c}%
y\leq-3\\
-3<y\leq-2\\
-2<y\leq1\\
1<y\leq3\\
3<y
\end{array}
.
\end{equation}
Some of the intermediary steps are illustrated in Figure \ref{HWfigure}.

Lemma 1 says that for a given $\left(  x,t\right)  $, $F\left(  y;x,t\right)
$ has a finite lower bound. \ This is trivial: for a given point $\left(
x,t\right)  $, as the generalised potential takes on only (at most) five
distinct nonzero values. \ Next we want to define $v\left(  x,t\right)  $ and
the set $S\left(  x,t\right)  .$ \ We have%
\begin{equation}
v\left(  x,t\right)  =\min_{y}F\left(  y;x,t\right)  =\min\left\{  \alpha
_{i}\right\}  _{i=1}^{5}%
\end{equation}
where $\alpha_{i}$ are the values that the piecewise function $F$ may take on,
and define%
\begin{equation}
S\left(  x,t\right)  =\left\{  y|\text{there exists }y_{n}\rightarrow y\text{
s.t. }F\left(  y_{n};x,t\right)  \rightarrow v\left(  x,t\right)  \right\}  .
\end{equation}
In this example, $v\left(  0,1\right)  =-\frac{1}{2}$ and for any $y\leq
x_{1}^{0}=-3$, we can clearly find $y_{n}\rightarrow y$ from the left so that
this condition holds. \ In fact, $F\left(  y;0,1\right)  $ achieves its
minimum at any such $y$, so that
\begin{equation}
v\left(  0,1\right)  =F\left(  y;0,1\right)  =-\frac{1}{2}%
\end{equation}
.

Lemma 2 says that for a point $y_{0}\in S\left(  x,t\right)  $ ($y_{0}\leq-3$
in our example) and $\left[  m_{0}\left(  y_{0}\right)  \right]  =m\left(
y_{0}+0\right)  -m\left(  y_{0}-0\right)  >0$ (which forces $y_{0}=-3$), then%
\begin{align}
v\left(  x,t\right)   & =\min_{y}F\left(  y;x,t\right)  =\left\{
\begin{array}
[c]{c}%
F\left(  y_{0};x,t\right) \\
F\left(  y_{0}+0;x,t\right)
\end{array}
\right.
\begin{array}
[c]{c}%
if\text{ }x\leq y_{0}+tu_{0}\left(  y_{0}\right) \\
if\text{ }x>y_{0}+tu_{0}\left(  y_{0}\right)
\end{array}
\nonumber\\
v\left(  0,1\right)   & =\min_{y}F\left(  y;0,1\right)  =F\left(
y_{0}+0;0,1\right)  =-\frac{1}{4}%
\end{align}
for $\left(  x,t\right)  =\left(  0,1\right)  ,$ since $0>-3+2.$

Lemma 3 says that for $\left(  x_{n},t_{n}\right)  $ and $y_{n}\in S\left(
x_{n},t_{n}\right)  $ that converge to $\left(  x,t\right)  $ and $y_{0},$
respectively, then $y_{0}\in S\left(  x,t\right)  $. \ In our example, this
corresponds to choosing $\left(  x_{n},t_{n}\right)  \rightarrow\left(
0,1\right)  $ and $y_{n}\in S\left(  x_{n},t_{n}\right)  .$ \ Recall%
\begin{equation}
F\left(  y;x_{n},t_{n}\right)  =\left\{
\begin{array}
[c]{c}%
\frac{3}{4}t_{n}+\frac{1}{4}\left(  x_{1}^{0}-x_{n}\right)  +\frac{1}%
{4}\left(  x_{2}^{0}-x_{n}\right) \\
\frac{1}{4}t_{n}+\frac{1}{4}\left(  x_{2}^{0}-x_{n}\right) \\
0\\
-\frac{1}{6}t_{n}+\frac{1}{3}\left(  x_{3}^{0}-x_{n}\right) \\
\frac{1}{3}\left(  x_{3}^{0}-x_{n}\right)  +\frac{1}{6}\left(  x_{4}^{0}%
-x_{n}\right)
\end{array}
\right.
\begin{array}
[c]{c}%
y\leq-3\\
-3<y\leq-2\\
-2<y\leq1\\
1<y\leq3\\
3<y
\end{array}
\end{equation}
Here we will have $S\left(  x_{n},t_{n}\right)  =\left\{  y\leq-3\right\}
=S\left(  x,t\right)  $ for $\left(  x_{n},t_{n}\right)  $ sufficiently close
to $\left(  x,t\right)  =\left(  0,1\right)  $. \ Indeed,%
\begin{equation}
F\left(  y;\varepsilon,1+\delta\right)  =\left\{
\begin{array}
[c]{c}%
-\frac{1}{4}+\frac{3}{4}\delta-\frac{1}{2}\varepsilon\\
-\frac{1}{8}+\frac{1}{4}\delta-\frac{1}{4}\varepsilon\\
0\\
\frac{1}{12}-\frac{1}{6}\delta-\frac{1}{3}\varepsilon\\
\frac{5}{12}-\frac{1}{2}\varepsilon
\end{array}
\right.
\begin{array}
[c]{c}%
y\leq-3\\
-3<y\leq-2\\
-2<y\leq1\\
1<y\leq3\\
3<y
\end{array}
\end{equation}
\ So clearly if $\left\{  y_{n}\right\}  $ converges, its limit is also in
$S\left(  x,t\right)  $.

\section{Introduction of Random Initial Conditions and Generalisation to
Conservation Laws}

Thus far, we have seen several examples of the evolution of conservation laws
under purely deterministic initial conditions. Other questions of interest
concern the behavior of the solutions under random initial conditions. It is
typical to assign an initial condition in the form of a stochastic process
with some structure, for example a Markov (\cite{A}, p. 144) or Feller
process. A Feller process (\cite{A}, p. 150) is essentially a special kind of
Markov process whose transition kernel is built from a semigroup with the
contraction property. A natural first question arising from this is whether
there are universality classes for which the structure of the initial
condition is preserved under the conservation law as time evolves \cite{MP}.
It is also of interest whether one can expect certain properties to persist in
time. For example, the Markov property in the continuum states that
information about the solution at one point offers no additional insight into
the solution at another, and is a very desirable condition to work with in the
context of probability theory. However, even if the Markov property is imposed
on the initial conditions, it may not persist for even an infinitesimal amount
of time. Another direction of interest is the extension of these results to
more general $C^{1}$, convex flux functions, beyond the special case of
Burgers' equation, where $f\left(  u\right)  =u^{2}/2$. Among many excellent
references for a clear explanation of several kinds of these stochastic
processes are \cite{A, BE, SC}.

Several key results along these lines are developed in \cite{M, MS}. First,
results are obtained involving the Lax equation from multiple different
perspectives. They note that a stationary, spectrally negative Feller process
can be characterised by a generator $A\left(  t\right)  $ acting on smooth
test functions $\varphi\in C_{c}^{1}\left(  \mathbb{R}\right)  $:%
\begin{equation}
A\varphi\left(  y\right)  =b\left(  y,t\right)  \varphi^{\prime}\left(
y\right)  +\int_{-\infty}^{y}\left(  \varphi\left(  z\right)  -\varphi\left(
y\right)  \right)  n\left(  y,dz,t\right)
\end{equation}
where $b\left(  y,t\right)  $ is the drift term and $n\left(  y,dz,t\right)  $
describes the jump density and satisfies%
\begin{equation}
\int1\wedge|y-z|^{2}n\left(  y,dz,t\right)  <\infty\text{ for all }%
y\in\mathbb{R}.
\end{equation}
By introducing a second operator which involves the flux function $f$, defined
by%
\begin{equation}
B\varphi\left(  y\right)  =-f^{\prime}\left(  y\right)  b\left(  y,t\right)
\varphi^{\prime}\left(  y\right)  -\int_{-\infty}^{y}\frac{f\left(  y\right)
-f\left(  z\right)  }{y-z}\left(  \varphi\left(  z\right)  -\varphi\left(
y\right)  \right)  n\left(  y,dz,t\right)  ,
\end{equation}
they show that evolution of the process is governed by the Lax equation%
\begin{equation}
\partial_{t}A=\left[  A,B\right]  =AB-BA.\label{Laxeqn}%
\end{equation}
These equations hold not only for the special case of Burgers' equation, as
first shown in \cite{CD}, but for more general fluxes. The structure of the
result is shown to depend only on the convexity assumption on the flux
function $f$ rather than Burgers' equation being a special case.

Another important result expands on kinetic theory. By expanding out the
commutator in (\ref{Laxeqn}), equations that describe shock clustering are
obtained. Namely, the drift satisfies%
\begin{equation}
\partial_{t}b\left(  y,t\right)  =-f^{\prime\prime}\left(  y\right)
b^{2}\left(  y,t\right)  ,
\end{equation}
and the jump density%
\begin{align}
& \partial_{t}n\left(  x,y,t\right)  +\partial_{y}\left(  nV_{y}\left(
y,z,t\right)  \right)  +\partial_{z}\left(  nV_{z}\left(  y,z,t\right)
\right) \nonumber\\
& =Q\left(  n,n\right)  +n\left(  \left(  \frac{f\left(  y\right)  -f\left(
z\right)  }{y-z}-f^{\prime}\left(  y\right)  \right)  \partial_{y}%
b-bf^{\prime\prime}\left(  y\right)  \right) \label{kteqn}%
\end{align}
where the velocities $V_{y}$ and $V_{z}$ are prescribed by%
\begin{equation}
V_{y}\left(  y,z,t\right)  =\left(  \frac{f\left(  y\right)  -f\left(
z\right)  }{y-z}-f^{\prime}\left(  y\right)  \right)  b\left(  y,t\right)
,\text{ }V_{z}\left(  y,z,t\right)  =\left(  \frac{f\left(  y\right)
-f\left(  z\right)  }{y-z}-f^{\prime}\left(  z\right)  \right)  b\left(
z,t\right)
\end{equation}
and $Q$ is a collision kernel describing various interactions between shocks.

It is proven that the equation (\ref{kteqn}) and the Lax equation
(\ref{Laxeqn}) are in fact equivalent. They also prove a closure property,
stating that if initial conditions are strong Markov (that is, have the Markov
property with respect to stopping times, \cite{A}, p. 97) and satisfy some
other additional assumptions, then a Markov property persists in time. More
precisely, one has (Thm 2, \cite{MS}):

\begin{theorem}
Define the inverse Lagrangian process $a\left(  x,t\right)  $ by the
following:%
\begin{align}
I(s;x,t)  & =\int_{0}^{s}\left(  u_{0}\left(  r\right)  -\left(  f^{\prime
}\right)  ^{-1}\left(  \frac{x-r}{t}\right)  \right)  dr,\nonumber\\
a\left(  x,t\right)   & =\arg^{+}\min_{s\in\mathbb{R}}I\left(  s;x,t\right)  ,
\end{align}
and let $u_{0}$ be a spectrally negative strong Markov process such that the
growth condition%
\begin{equation}
\lim_{\left\vert s\right\vert \rightarrow\infty}I\left(  s;x,t\right)
=+\infty
\end{equation}
holds a.s. Under the law $%
\mu
_{0}$ and for any fixed $t>0$, the inverse Lagrangian process $a\left(
x,t\right)  $ is Markov.
\end{theorem}

\section{Exact Solution For a Special Case of Random Initial Conditions}

With the introduction of random initial conditions into Burgers' equation and
more general conservation laws, in general, the most one can hope to obtain
for a solution in terms of complicated generators, as we have seen in Section
5. However, for a few specialised kinds of initial data, one can obtain exact,
closed-form solutions up to the level of special functions. One such case is
that of Frachebourgh and Martin \cite{FM}, for which white noise initial data
is considered, and expressions for the one- and two-point functions obtained,
building off a key result of \cite{GR}. The closed analytical forms are
expressed in terms of Airy functions. The starting point involves the inviscid
Burgers' equation%
\begin{equation}
\frac{\partial}{\partial t}u\left(  x,t\right)  +u\left(  x,t\right)
\frac{\partial}{\partial x}u\left(  x,t\right)  =\nu\frac{\partial}{\partial
x^{2}}u\left(  x,t\right) \label{viscidburg}%
\end{equation}
and using variational methods in tandem with taking the limit of viscosity
parameter $\nu\downarrow0$. To this end, one may introduce the potential
$\partial\Psi\left(  x,t\right)  /\partial x=u\left(  x,t\right)  $ and use
the Cole-Hopf transformation $\Psi\left(  x,t\right)  =-2\nu\ln\theta\left(
x,t\right)  $ along with other methods (\cite{EV}, p. 207) to show that
$\theta\left(  x,t\right)  $ satisfies the heat equation. Consequently, the
solution to (\ref{viscidburg}) is then given by%
\begin{equation}
u\left(  x,t\right)  =\frac{\int_{-\infty}^{\infty}dy\frac{x-y}{t}\exp\left(
-\frac{1}{2\nu}F\left(  x,y,t\right)  \right)  }{\int_{-\infty}^{\infty}%
dy\exp\left(  -\frac{1}{2\nu}F\left(  x,y,t\right)  \right)  }%
\label{viscidsol}%
\end{equation}
where%
\begin{align}
F\left(  x,y,t\right)   & =\frac{\left(  x-y\right)  ^{2}}{2t}-\psi\left(
y\right) \nonumber\\
\psi\left(  y\right)   & =-\Psi\left(  y,0\right)  =-\int_{0}^{y}dy^{\prime
}u\left(  y^{\prime},0\right)  ,
\end{align}
the latter of which depends, obviously, on the given initial conditions. In
the limit of interest $\nu\downarrow0,$ the only contributions from
(\ref{viscidsol}) arise from where $F$ has a minimum, i.e. at point(s)
prescribed by%
\begin{equation}
\xi\left(  x,t\right)  =\min_{y}F\left(  x,y,t\right)  ,
\end{equation}
and write, formally%
\begin{equation}
u\left(  x,t\right)  =\frac{x-\xi\left(  x,t\right)  }{t}.
\end{equation}
Due to the scaling properties of the solution, it is sufficient to set $t=1$
herein. To find this minimum, $\xi\left(  x,1\right)  ,$ consider the
following geometric illustration of the solution. Picture $\psi\left(
x\right)  $ as a sample Brownian motion path, and a parabola $P\left(
x\right)  =\left(  x-y\right)  ^{2}/2+C$. One then adjusts $P\left(  x\right)
(=P_{y}\left(  x\right)  )$ by sliding it down, decreasing $y$ so that it
touches the Brownian path but does not cross it. This may occur at multiple
points, or at just one point. If there are two such points for a value
$x^{\ast}$, labelled $\xi_{-}$ and $\xi_{+},$ then the function $F\left(
x,y,t\right)  $ has a discontinuity at $x^{\ast}$, $\mu$ which is called a
\textit{shock}. We then characterise this shock by two parameters:%
\begin{equation}
\mu=\xi_{+}-\xi_{-}\text{, }\nu=x^{\ast}-\xi_{-},
\end{equation}
named \textit{strength} and \textit{wavelength}.
\begin{figure}
[h]
\begin{center}
\includegraphics[width=\linewidth]{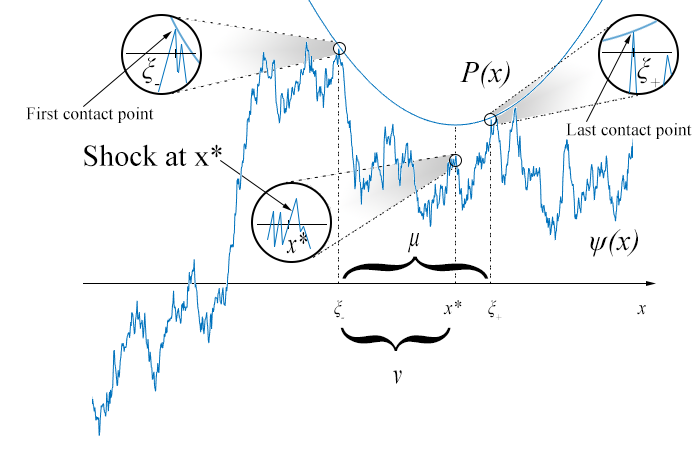}
\caption{For a shock at a point $x^{\ast}$, we slide the parabola $\left(
x-x^{\ast}\right)  ^{2}/2$ down until we have (at least) two contact points
with the Brownian path, but in such a way that the parabola does not cross the
Brownian path. If there are more than two, we consider only the first and last
contact points. These points are given by $\left(  \xi_{-},\left(  \xi
_{-}-x^{\ast}\right)  ^{2}/2\right)  $ and $\left(  \xi_{+},\left(  \xi
_{+}-x^{\ast}\right)  ^{2}/2\right)  $. The shock is then described by the
parameters $\mu=\xi_{+}-\xi_{-}$ and $\nu=x^{\ast}-\xi_{-}.$ This figure is
based off Figure 1, \cite{FM}.}
\label{FMfigure}
\end{center}
\end{figure}

These notions are illustrated in Figure \ref{FMfigure}. This forms the basis
of relating the one-point functions $p_{1}\left(  x,u\right)  ,$ the
probability that the velocity field at a point $x$ takes values between $u$
and $u+du,$ with the object $\rho_{1}\left(  \mu,\eta\right)  $. In
particular, one can adjust the coordinate system and has%
\begin{equation}
\rho_{1}\left(  \mu,\eta\right)  =\mathbb{E}\left\{
\begin{array}
[c]{c}%
\psi\left(  x\right)  \leq P_{\nu}\left(  x\right)  ,x\in\mathbb{R},\text{
first contact with }P_{\nu}\left(  x\right)  \text{ at }\left(  0,0\right)
;\\
\text{last contact with }P_{\nu}\left(  x\right)  \text{ at }\left(  \mu
,\eta\right)
\end{array}
\right\}
\end{equation}
Similarly, one can consider two such parabolas and build the two-point
functions $p_{2}\left(  x_{1},u_{1},x_{2},u_{2}\right)  $, giving the
probability that the velocity field at points $x_{1}$, $x_{2}$ take values
between $u_{1}$ and $u_{1}+du$, $u_{2}$ and $u_{2}+du_{2}$, in terms of
$\rho_{2}\left(  0,\mu_{1},\eta_{1},x,\mu_{2},\eta_{2}\right)  $.

Integrating up the shock strength distribution $\rho_{1}\left(  \mu\right)
=\int_{-\infty}^{\infty}\rho_{1}\left(  \mu,\eta\right)  d\eta,$ they obtain
that the shock distribution as a function of strength is given by%
\begin{equation}
\rho_{1}\left(  \mu\right)  =2a^{3}\mu\sum_{k\geq1}e^{-a\omega_{k}\mu}\frac
{1}{\left(  2\pi\right)  ^{2}}\int_{-\infty}^{\infty}d\zeta_{1}\int_{-\infty
}^{\infty}d\zeta_{2}\frac{e^{-ia\mu\left(  \zeta_{1}+\zeta_{2}\right)  /2}%
}{Ai\left(  i\zeta_{1}\right)  Ai\left(  i\zeta_{2}\right)  }\int_{-\infty
}^{\infty}d\eta^{\prime}e^{i\eta^{\prime}\left(  \zeta_{1}-\zeta_{2}\right)  }%
\end{equation}
In the same vein, they obtain results for $\rho_{2}$, thus giving an exact,
closed-form expression for the characterisation of the system dynamics in this
case of Burgers' equation with white noise initial data.

\section{Further Extension of Results to Flux Functions With Less Regularity}

As we have seen in the previous sections, conservation laws such as Burgers'
equation can admit an exact solution in some special cases, even with random
initial data, as in \cite{FM}. More general expressions entailing generators
and semigroup theory can be extended to the case of an arbitrary $C^{1}$ flux
function in \cite{MS} and provide more insight into the solution. Using a
specific set of test functions and probability theory, we show that one may
consider a nonlinear flux function that is only continuous and obtain
meaningful results for the hierarchy of equations describing the formation and
interaction of shocks. The results obtain are in terms of various n-point
functions, representations of probabilities (or probability densities) that
the solution field takes certain values at a given time and a number of
different positions.

\subsection{Density of States Approach}

\begin{figure}
[ptb]
\begin{center}
\includegraphics[width=\linewidth]{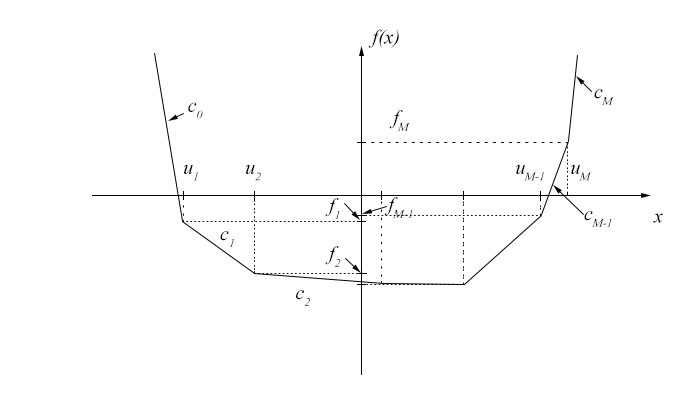}
\caption{Illustration of the flux function as described above.}
\label{Caginalp1}
\end{center}
\end{figure}

We build the solution to a piecewise linear flux function defined by piecewise
linear interpolation between the points%
\begin{equation}
f\left(  u_{i}\right)  =f_{i},\text{ }1\leq i\leq M
\end{equation}
with slopes consequently given by%
\begin{equation}
c_{k}=\frac{f_{k+1}-f_{k}}{u_{k+1}-u_{k}}.
\end{equation}
The flux function is illustrated in Figure \ref{Caginalp1}. Burgers' equation%
\begin{equation}
u_{t}+uu_{x}=0\label{burgers}%
\end{equation}
can be written in its entropy-entropy flux pair form in terms of smooth test
functions $\varphi$ and $\psi$ as%
\begin{equation}
\partial_{t}\varphi\left(  u\left(  x,t\right)  \right)  =-\partial_{x}%
\psi\left(  u\left(  x,t\right)  \right)  .\label{1ptstart}%
\end{equation}

Further technical details of the equivalence between the entropy-entropy flux
pair form (\ref{1ptstart}) and (\ref{burgers}) can be found in \cite{VO}.

Taking expectations of (\ref{1ptstart}) and formally passing derivatives
inside yields%
\begin{equation}
\mathbb{E}\left\{  \partial_{t}\varphi\left(  u\left(  x,t\right)  \right)
\right\}  =\partial_{t}\mathbb{E}\left\{  \varphi\left(  u\left(  x,t\right)
\right)  \right\}  =-\mathbb{E}\left\{  \partial_{x}\psi\left(  u\left(
x,t\right)  \right)  \right\}  .\label{1ptexp}%
\end{equation}
As this is a piecewise linear function with $M-1$ distinct slopes, it makes
sense to consider a test function $\varphi$ with discrete values on the set of
points $x\in\left\{  x_{i}\right\}  _{i=1}^{M},$ extended for other $x$ to
make it left continuous. One then has%
\begin{equation}
\partial_{t}\sum_{l=1}^{M}\varphi\left(  u_{l}\right)  p_{1}\left(
x,t;u_{l}\right)  =-\sum_{l=1}^{M}\sum_{m=1}^{M}\left(  \psi\left(
u_{m}\right)  -\psi\left(  u_{l}\right)  \right)  p_{2}\left(  x,x+,t;u_{l}%
,u_{m}\right)  .\label{1ptsum}%
\end{equation}

\begin{figure}
[ptb]
\begin{center}
\includegraphics[width=\linewidth]{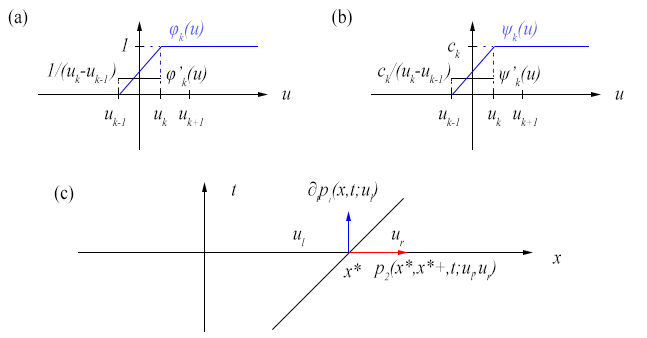}
\caption{(a)-(b)\ Construction of the test functions $\varphi_{k}\left(
u\right)  $ and $\psi_{k}\left(  u\right)  $; (c) Illustration of a shock,
with positive contribution from $\partial_{t}p_{1}\left(  x,t;u_{l}\right)  $
(upward arrow, blue), and negative contribution from $p_{2}\left(
x,x+,t;u_{l},u_{r}\right)  $ (right arrow, red).}
\label{Caginalp2}
\end{center}
\end{figure}

We now choose the derivative of the test function $\varphi$ to be a
discretised version of a $\delta$-function as $u_{k}$ for a given $k$, i.e.%
\begin{equation}
\varphi_{k}^{\prime}\left(  u\right)  =\frac{1_{[u_{k,},u_{k+1})}\left(
u\right)  }{u_{k+1}-u_{k}},
\end{equation}%
\begin{equation}
\varphi_{k}\left(  u\right)  =1_{[u_{k+1},\infty)}\left(  u\right)  ,\text{
}\psi_{k}^{\prime}\left(  u\right)  =\frac{c_{k}1_{[u_{k},u_{k+1})}\left(
u\right)  }{u_{k+1}-u_{k}},\text{ }\psi\left(  u\right)  =c_{k}1_{[u_{k+1}%
,\infty)}\left(  u\right)  .\label{1ptint}%
\end{equation}
These are further illustrated in Figure \ref{Caginalp2}(a-b). Substituting the
expressions (\ref{1ptint}) into (\ref{1ptsum}), we obtain%
\begin{equation}
\sum_{l=k+1}^{M}\partial_{t}p_{1}\left(  x,t;u_{l}\right)  =\sum_{l=k+1}%
^{M}\sum_{m=1}^{k}c_{k}p_{2}\left(  x,x+,t;u_{l},u_{m}\right)  -\sum_{l=1}%
^{k}\sum_{m=k+1}^{M}c_{k}p_{2}\left(  x,x+,t;u_{l},u_{m}\right)
\label{1ptfirstmethodfinal}%
\end{equation}

This is significant as it provides a relation for the change (in time)\ of the
value of the one-point distribution at position $x$\ as a function of the
difference of two terms involving the two-point distribution at the same point
$x$ and the value of the solution at its right limit $x+$. This difference of
terms can be interpreted as a discrete derivative or a finite difference. By
generalising the procedure outlined above, we obtain similar results for the
hierarchy at the n-point level. Thus, the expression for the n-point function
at points $x_{1},...,x_{n}$ depends on analogous pairs of terms with interplay
between value of the $n+1$ point function at $x_{i}+$ and $x_{i}$. The
relation (\ref{1ptfirstmethodfinal}) applies uniformly before shock
interactions, but the meaningful case one should visualise is when a shock is
present at position $x$. This is shown in Figure \ref{Caginalp2}(c).

\subsection{Density of Shocks Approach}

Although the first approach outlined above accurately describes the initial
formation and propagation of shocks, it requires additional techniques to
persist through shock collisions. By considering a second method and changing
the interpretation of the n-point function to be inclusive of both the value
of the velocity field $u$ at $x$ and $x+,$ one is able to derive a second
hierarchy that accurately described the system dynamics even through these
interactions. Specifically, the new definition linked the information of the
value of the solution at both limits, making $x_{i}$ and $x_{i}+$ always
appear together, replacing the old expressions as follows:%
\begin{align}
p_{1}\left(  x,t;u\right)   &  \rightarrow f_{1}\left(  x_{-},x,t;u,v\right)
\nonumber\\
p_{2}\left(  x_{1},x_{2},t;u_{1},u_{2}\right)   &  \rightarrow f_{2}\left(
x\,_{1,-},x_{1},x_{2,-},x_{2},t;u_{1},v_{1},u_{2},v_{2}\right)
.\label{new1pt}%
\end{align}

In deriving this hierarchy, it clearly makes sense to have a
\textit{transport} term of the form
\begin{equation}
\partial_{t}f_{1}\left(  x_{-},x,t;u,v\right) \label{1pttimederiv}%
\end{equation}
along with a \textit{free-streaming} term%
\begin{equation}
c_{uv}\partial_{x}f_{1}\left(  x_{-},x,t;u,v\right)  .\label{1ptstream}%
\end{equation}
When considering the hierarchy at the nth level, (\ref{1ptstream}) is
generalised by having a sum of $n$ terms, one for each shock between $u_{i}$
and $v_{i}$. On the other side of equation, terms for various interactions
causing creation or destruction of the shock $u$ and $v$ are added. Using
Taylor series, one picks up appropriate rate constants and obtains%
\begin{align}
&  \partial_{t}f_{1}\left(  x,t;u,u+1\right)  +c_{u,u+1}\partial_{x}%
f_{1}\left(  x,t;u,u+1\right) \nonumber\\
&  =-\sum_{w<u}\left(  c_{u,u+1}-c_{u+1,w}\right)  \partial_{2}f_{2}\left(
x,x,t;u,u+1,u+1,w\right) \nonumber\\
&  -\sum_{w>u+1}\left(  c_{w,u}-c_{u,u+1}\right)  \partial_{1}f_{2}\left(
x,x,t;w,u,u,u+1\right) \label{2nd1pt1}%
\end{align}
for $u<v$ (for which $v=u+1$ yields the only nontrivial case) and
\begin{align}
&  \partial_{t}f_{1}\left(  x,t;u,v\right)  +c_{uv}\partial_{x}f_{1}\left(
x,t;u,v\right) \nonumber\\
&  =\sum_{w\leq u+1}\left(  c_{uw}-c_{wv}\right)  \partial_{2}f_{2}\left(
x,x,t;u,w,w,v\right)  -\sum_{w\leq v+1}\left(  c_{uv}-c_{vw}\right)
\partial_{2}f_{2}\left(  x,x,t;u,v,v,w\right) \nonumber\\
&  -\sum_{w\geq u-1}\left(  c_{wu}-c_{uv}\right)  \partial_{1}f_{2}\left(
x,x,t;w,u,u,v\right)  .\label{2nd1pt2}%
\end{align}
for $v<u$. We use the notation $u+1$ to mean the nearest neighbor state, with
the states in ascending order. When tested against examples, it is confirmed
that these equations accurately model the dynamics after shocks collide
without modifications, unlike the first approach or in many classical results.

A related problem is studied in \cite{KR}, where an expression for the
statistics is obtained in the problem with an initial condition bounded to the
case $0\leq x\leq L$. In an approach reminiscent of statistical mechanics
methods, limits of a large number of particles are taken to obtain these results.

\section{Open Problems and Ideas For Future Research}

There remain a number of open problems, particularly with the introduction of
randomness into the initial conditions for conservation laws. Additional
problems include the more detailed analysis and perhaps direct computation
with broader classes of initial conditions. Although the Markov property
imposed on random initial conditions does not persist under the evolution of
time for Burgers' equation, other classes of random initial data and/or more
general conservation laws can be considered and explored further.

A closure theorem for the n-point function is another possible direction for
future work. Such a result would allow the time derivative of the n-point
function to be expressed explicitly in terms of the spatial derivatives of
n-point functions of different arguments. It is possible, for example, that
the n+1-point functions may factor into lower point functions for larger
values of $n$. Such a reduction, should it exist, would serve to greatly
enhance the understanding of conservation laws from yet another perspective.

\bibliographystyle{plain}
\bibliography{mybib}

\end{document}